\documentclass[10pt]{article}
\usepackage{graphicx}
\usepackage{amssymb}
\usepackage{epstopdf}
\usepackage{amsmath}
\usepackage{amsfonts}
\newcommand{\braket}[2]{\langle #1,#2 \rangle}

\def\phi{{\varphi}}

\DeclareSymbolFont{AMSb}{U}{msb}{m}{n}
\DeclareMathSymbol{\N}{\mathbin}{AMSb}{"4E}
\DeclareMathSymbol{\Z}{\mathbin}{AMSb}{"5A}
\DeclareMathSymbol{\R}{\mathbin}{AMSb}{"52}
\DeclareMathSymbol{\Q}{\mathbin}{AMSb}{"51}
\DeclareMathSymbol{\I}{\mathbin}{AMSb}{"49}
\DeclareMathSymbol{\C}{\mathbin}{AMSb}{"43}

\begin{document}

\addtolength{\textheight}{0 cm} \addtolength{\hoffset}{0 cm}
\addtolength{\textwidth}{0 cm} \addtolength{\voffset}{0 cm}

\newenvironment{acknowledgement}{\noindent\textbf{Acknowledgement.}\em}{}

\setcounter{secnumdepth}{5}
 \newtheorem{proposition}{Proposition}[section]
\newtheorem{theorem}{Theorem}[section]
\newtheorem{lemma}[theorem]{Lemma}
\newtheorem{coro}[theorem]{Corollary}
\newtheorem{remark}[theorem]{Remark}
\newtheorem{ex}[theorem]{Example}
\newtheorem{claim}[theorem]{Claim}
\newtheorem{conj}[theorem]{Conjecture}
\newtheorem{definition}[theorem]{Definition}
\newtheorem{application}{Application}

\newtheorem{corollary}[theorem]{Corollary}

\title{A Self-dual Polar Factorization for Vector Fields}
\author{Nassif  Ghoussoub\thanks{Partially supported by a grant
from the Natural Sciences and Engineering Research Council of Canada.}\\
{\it\small Department of Mathematics}\\
{\it\small  University of British Columbia}\\
{\it\small Vancouver BC Canada V6T 1Z2}\\
{\it\small nassif@math.ubc.ca}\vspace{1mm}\and
Abbas Moameni\thanks{Research supported by a Coleman fellowship at Queen's  University.}\hspace{2mm}\\
{\it\small Department of Mathematics and Statistics}\\
{\it\small Queen's University}\\
{\it\small Kingston, ON,  Canada K7L3N6}\\
{\it\small  momeni@mast.queensu.ca}\\
\date{Revised August 1, 2011}
%\today
}

\maketitle

\begin{abstract} We show that any non-degenerate vector field $u$ in $ L^{\infty}(\Omega, \R^N)$, where $\Omega$ is a bounded domain in $\R^N$, can be written as
\begin{equation}
\hbox{$u(x)= \nabla_1 H(S(x), x)$ for a.e. $x \in \Omega$},
\end{equation}
where $S$ is  a measure preserving point transformation on $\Omega$ such that $S^2=I$ a.e (an involution), and $H: \R^N \times \R^N \to \R$ is a globally Lipschitz  anti-symmetric convex-concave Hamiltonian. Moreover, $u$ is a monotone map if and only if $S$ can be taken to be the identity, which suggests that our result is a self-dual version of Brenier's polar decomposition for the vector field $u$ as  $u(x)=\nabla \phi (S(x))$, where $\phi$ is convex and $S$ is a measure preserving transformation. We also describe how our polar decomposition can be reformulated as a self-dual mass transport problem.
\end{abstract}

\section{Introduction}

Given a bounded domain $\Omega$ in $\R^N$, a classical theorem of Rockafellar \cite{Ph} yields that
 a single-valued map $u$ from $\Omega$ to $\R^n$  is  {\it cyclically monotone}, i.e.,  for any finite number of points $(x_i)_{i=0}^n$  in $\Omega$ with $x_0=x_n$, we have
\begin{equation}
\hbox{$\sum\limits_{i=1}^n\langle u(x_k), x_k-x_{k-1} \rangle \geq 0$,}
\end{equation}
if and only
\begin{equation}
\hbox{$u (x)=\nabla \phi (x)$ on $\Omega$, }
\end{equation}
where $\phi:\R^n \to \R$ is a convex function.

On the other hand, a result of E. Krauss  \cite{Kra} yields that  $u$ is a {\it   monotone map}, that is
%$u+\lambda I$ is invertible and
for all $(x,y)$ in $\Omega$, \begin{equation}
\hbox{$\langle x-y, u(x) - u(y)\rangle \geq 0$,}
\end{equation}
  if and only if
  \begin{equation}
  \hbox{$u(x)=\nabla_1H(x, x)$ for all $x\in \Omega$,}
  \end{equation}
   where $H$ is a  convex-concave anti-symmetric Hamiltonian on $R^N\times R^N$.

More remarkable is the polar decomposition that Y. Brenier \cite{Br} establishes for a general  non-degenerate vector field, and which follows from his celebrated mass transport theorem. Recall that a mapping $u: \Omega \rightarrow \R^N$ is said to be {\it non-degenerate} if  the inverse image $u^{-1}(N)$ of  every zero-measure $N\subseteq \R^N$ has also zero measure.

Brenier then proved stating that any non-degenerate vector field $u \in L^{\infty}(\Omega, \R^N)$ can be decomposed  as
\begin{equation}
\hbox{$u(x)=\nabla \psi \circ S(x)$ \,\, a.e. in $\Omega$,}
\end{equation}
with
$\psi:\R^N\rightarrow \R$  being a  convex function and  $S:\bar {\Omega}\rightarrow \bar {\Omega}$ is a measure preserving
transformation.

In this paper, we shall prove another decomposition for non-degenerate vector fields, in  the same spirit of Brenier's, but which can be seen as the general version of Krause's. Indeed, here is the  main result of this paper.

\begin{theorem}\label{main} Let $\Omega$ be an open bounded set in $\R^N$ such that $meas(\partial \Omega)=0$.
\begin{enumerate}
\item  If $u \in L^{\infty}(\Omega, \R^N)$ is a non-degenerate vector field,  then there exists a measure preserving transformation $S:\bar
{\Omega}\rightarrow \bar {\Omega}$ such that $S^2=I$ (i.e., an involution), and a globally Lipschitz  anti-symmetric convex-concave Hamiltonian $H: \R^N \times \R^N \to \R$ such that
\begin{equation} \label{rep.1}
u(x)= \nabla_1 H(Sx, x) \quad \text{ a.e. } \, x \in \Omega.
\end{equation}
\item If $u$ is differentiable a.e., and the map
\begin{equation}\label{unique}
 x\to \langle u(x), y_1-y_2\rangle +\langle u(y_1)-u(y_2), x\rangle
\end{equation}
has no critical points in $\Omega$ unless $y_1=y_2$, then there exists a unique measure preserving involution $S$ such that satisfies $(\ref{rep.2})$ for some convex-concave anti-symmetric Hamiltonian $H$.
\item In particular,  if $u$ is a strictly monotone map, then $S$ is necessarily equal to the idendity.
\end{enumerate}
\end{theorem}
Since $S$ is an involution and $H$ is anti-symmetric, one can deduce from (\ref{main}) that
\begin{equation} \label{Ham.2}
u(Sx)= -\nabla_2 H(Sx, x) \quad \text{ a.e. } \, x \in \Omega.
\end{equation}

The fact that $S$ is a measure preserving involution provides an improvement on the first factor in Brenier's decomposition. On the other hand,  the second factor i.e., the potential $\nabla \psi$, is obviously better than  the partial gradient of a convex-concave anti-symmetric Hamiltonian on $\R^N \times \R^N$, which is only a 
monotone map.
The connection to self-duality will be described later in this introduction.

We now give a few examples of how this decomposition appears in concrete situations.\\

\textbf{1. Basic monotone operators not derived from a potential:} If $u=\nabla \phi + A$, where $\phi$ is convex and $A$ is a skew-adjoint matrix, then clearly
\begin{equation}
H(x,y)=\phi (x)-\phi (y) -\langle Ax, y\rangle,
\end{equation}
and $S$ is the identity. This is of course a very important case of maximal monotone operators that we shall discuss later. \\

\textbf{2. Helmholtz Decomposition of vector fields:} Let $\Omega$ be a smooth bounded connected open set in $\mathbb{R}^N.$ Any smooth vector field $u$ on $\Omega$ can then be written in a unique way as $u(x)=\nabla p(x)+ v(x),$ where $p$ is a smooth real function on $\Omega$, and $v$ is a smooth divergence free vector field parallel to the boundary of $\Omega$.
By considering $u_\epsilon$ as a smooth perturbation of the identity map:
\[u_{\epsilon}(x)=x +\epsilon u(x), \qquad x \in \bar \Omega.\]
for $\epsilon$ small enough, we can write
\[u_{\epsilon}(x) =\nabla_1 H(x,x), \qquad \text{ for all } x \in \bar \Omega,\]
where
\[H(x,y)= \frac{1}{2}|x|^2+\epsilon p(x)+ \epsilon \langle x-y, \frac{v(x)+v(y)}{2}\rangle - \frac{1}{2}|y|^2-\epsilon p(y),\]
and again $S$ is the identity operator. Note that for $\epsilon$ small enough, $H$ is convex in the first variable and concave in the second one.\\

Note that both examples above fit in the framework of the result of E. Krauss  \cite{Kra}, who --as mentioned above-- has shown that if $u$ is a single-valued maximal monotone map on $\mathbb{R}^N$, then $u(x)=\nabla_1 H(x,x)$, that is $S$ is the identity map whenever $u$ is a monotone map.  We shall come back to this theme when we discuss self-duality below. Now, we give examples of non-monotone operators. \\

\textbf{3. Decomposition of real matrices:} Any $N \times N$ matrix $A$ can be decomposed  as $A_s +A_a$ where $A_s$ is the symmetric part and $A_a$ is the anti-symmetric part. The symmetric part $A_s$ can then  be written as the product $RS$ of a symmetric non-negative matrix $R$ and  a real unitary matrix $S$. It is not difficult to check that since $A_s$ is symmetric so is $S$ and therefore $S^2=I$.

It follows that we can write the following decomposition for the matrix $A:$
\[Ax =A_ax+RSx=\nabla_1 H(S(x),x), \qquad \text{ for all } x \in R^N,\]
where
\[
H(x,y)=\frac{1}{2} \langle R x,x \rangle  -\frac{1}{2} \langle R y,y \rangle -\langle A_a x, y \rangle
\]
is clearly a skew-adjoint Hamiltonian and $S$ is a symmetric involution matrix.\\
Note that the condition (\ref{unique}) is, in this case, equivalent to saying that the symmetric part
$A_s$ of $A$ is non-singular. Indeed,
\[
\langle u(x), y_1-y_2\rangle +\langle u(y_1)-u(y_2), x\rangle =2\langle A_s(y_1-y_2), x\rangle
\]
has a critical point if only if $A_sy_1=A_sy_2$, which means that $y_1=y_2$ whenever $A_s$ is assumed to be non-singular. This is compatible with the classical fact that  $R$ and $S$ in the above decomposition of $A_s$  are unique. \\

\textbf{4. Examples of representations on non-monotone maps on the line:}\\ (i) A simple non-monotone example is the function \begin{equation}
u(x)=\sin x +x \cos x.
\end{equation}
 It can be written as $u(x)=\nabla_1 H(S(x),x)$ on $[0, \pi]$, where
\begin{equation}
\hbox{$H(x,y)= x\sin y-y\sin x$ \quad and  \quad $S(x)=\pi-x$.}
\end{equation}
(ii) More generally, a large class of examples of convex-concave anti-symmetric Hamiltonians is given by
\[
H(x,y)=f(x)g(y)-f(y)g(x) +h(x-y),
\]
where $h$ is odd and with suitable conditions that render $H$ convex in $x$. For $\alpha \in [0,1]$, the function
 \begin{eqnarray}
  \label{eleven}
u(x)=\left\{
  \begin{array}{ll}
 \hfill   f'(\alpha - x) g(x)-g'(\alpha -x)f(x)+h'(\alpha-2x)& \quad {\rm if}  \,\,\, 0\leq x\leq \alpha\\
 \hfill f'(x)g(x)-g'(x)f(x)+h'(0) & \quad {\rm if} \,\, \alpha < x\leq 1.  \end{array}
\right.
\end{eqnarray}
can then be written as $u(x)=\nabla_1 H(S(x),x)$, where $S(x)=\alpha-x$ on $[0, \alpha]$ and $S(x)=x$ on $(\alpha, 1)$. \\

(iii) A more interesting example is the map
  \begin{eqnarray}
  \label{eleven}
u(x)=\left\{
  \begin{array}{ll}
 \hfill   3-2x & \quad {\rm if}  \,\,\frac{1}{2}\leq x\leq 1\\
 \hfill 2x & \quad {\rm if} \,\, 0\leq x\leq \frac{1}{2}.  \end{array}
\right.
\end{eqnarray}
One can easily verify that $u(x)=\nabla_1 H(S_1(x),x)$, where $S_1(x)=1-x$ and $H$ is given by the following formula
 \begin{eqnarray}
  \label{eleven}
H(x,y)=\left\{
  \begin{array}{llll}
 \hfill   -2xy+2x-y-\frac{1}{2} & \quad {\rm if}  \,\, 0\leq x\leq \frac{1}{2}&{\rm and} & \frac{1}{2}\leq y\leq 1\\
 \hfill x-y & \quad {\rm if} \,\, 0\leq x, y\leq \frac{1}{2}\,\, &{\rm or} & \frac{1}{2}\leq x, y\leq 1 \\
  \hfill   2xy-2y+x+\frac{1}{2} & \quad {\rm if}  \,\,  \frac{1}{2}\leq x\leq 1 &{\rm and} & 0\leq y\leq \frac{1}{2}.
\end{array}
\right.
\end{eqnarray}
Also note that $u(x)=\nabla_1 H(S_2(x),x)$, where $S_2(x)=x+\frac{1}{2}$ on $[0, \frac{1}{2})$ and  $S_2(x)=x-\frac{1}{2}$ on $(\frac{1}{2}, 1]$, which means that the involution $S$ appearing in the decomposition is not necessarily unique.\\
Actually, one has non-uniqueness whenever there exists two subsets $A,B$ of positive measure such that
\[
\hbox{$ u(x) . y + u(y) . x = f(x) + g(y)$  when $ (x,y)\in  A \times B$,}
 \]
 for some functions $f$ and $g$. This is what happens in the previous example with $
A = [0, 1/2]$ and $B = [1/2, 1]$. It follows that $x\to u(x) . (y_2 - y_1) + (u(y_2) - u(y_1)) . x $ doesn't depend on $x$.  If $u$ is differentiable then $u'(x) + u'(y) = 0$  for $(x,y)\in A \times B$. \\

\textbf{5. Why can the decomposition be considered ``self-dual"?:} Let $X$ be a reflexive Banach space, and recall from  \cite{Gh}  the notion of a vector field $\bar \partial L$ that is derived from a convex lower semi-continuous Lagrangian on phase space $L:X\times X^*\to \R\cup\{+\infty\}$ in the following way:  for each $x\in X$, the  --possibly empty-- subset  $\bar \partial L(x)$ of $X^*$ is defined as
 \begin{eqnarray}
\bar \partial L(x): = \{ p \in X^*;  (p, x)\in \partial L(x,p) \}.
\end{eqnarray}
 Here $\partial L$ is the subdifferential of the convex function $L$ on $X\times X^*$, which should not be confused with $\bar \partial L$. Of particular interest are those vector fields derived from  {\it self-dual Lagrangians}, i.e., those convex lower semi-continuous Lagrangians $L$ on $X\times X^*$ that satisfy the following duality property:
  \begin{equation}
  L^*(p,x)=L(x, p) \quad \hbox{\rm for all $(x,p)\in X\times X^*$},
  \end{equation}
   where here $L^*$ is the Legendre transform in both variables, i.e.,
  \[
  L^*(p,x)= \sup  \{ \braket{ y}{p }+  \braket{x}{q }-L(y, q): \,  (y,q)\in X\times X^*\}.
  \]
  Such Lagrangians satisfy the following basic property:
 \begin{equation}\label{obs.1}
\hbox{$ L(x, p)-\langle x, p\rangle\geq 0$   for every $(x, p) \in X\times X^{*}$.}
 \end{equation}
  Moreover,
 \begin{equation}\label{obs.2}
\hbox{  $L(x, p)-\langle x, p\rangle =0$ if and only if $(p, x)\in \partial L(x,p)$,}
\end{equation}
which means  that the associated vector field at $x \in X $ is simply
\begin{eqnarray}
\bar \partial L(x):= \{ p \in X^*; L(x,p)- \langle x,p \rangle=0 \}.
\end{eqnarray}
These so-called  {\it selfdual vector fields} are natural but far reaching extensions of subdifferentials of convex lower semi-continuous functions. Indeed, the most
basic selfdual Lagrangians are of the form $L(x,p)= \varphi (x)+\varphi^*(p)$ where  $\varphi$ is such a function in $X$, and $\varphi^*$
is its Legendre conjugate on $X^*$,  in which case  $\bar \partial  L(x)= \partial \varphi (x)$.

More interesting  examples of self-dual Lagrangians are of the form $L(x,p)= \varphi (x)+\varphi^*(-\Gamma x+p)$  where $\varphi$  is a convex and lower semi-continuous
function on $X,$ and  $\Gamma: X\rightarrow X^*$ is a skew adjoint operator. The corresponding selfdual vector field is then
$\bar \partial L(x)=\Gamma x+ \partial \varphi (x)$.

 Actually, it turned out that any {\it maximal monotone operator} $A$ (a notion studied for example in  \cite{Br})  is a self-dual vector field and vice-versa  \cite{Gh}. That is,  there exists a selfdual Lagrangian $L$ such that $A=\bar \partial L$. This fact was proved and reproved by several authors. See for example, R.S. Burachik and B. F. Svaiter  \cite{BS}, B. F. Svaiter \cite{S}, and Baushke and Wang \cite{BW}).
 
This result means that self-dual Lagrangians can be seen as the {\it potentials} of maximal monotone operators, in the same way as the Dirichlet integral  is the potential of the Laplacian operator (and more generally as any convex lower semi-continuous energy is a potential for its own subdifferential). Check \cite{Gh} to see how this characterization leads to variational formulations and resolutions of most equations involving monotone operators.

 Consider now the  Hamiltonian $H_L$ on $X^*\times X^*$ corresponding to $L$, that is the Legendre transform of $L$ in the first variable,
 \[
 H_L(p,q)=\sup\{\langle x, p\rangle - L(x, q); x\in X\}.
 \]
It is convex-concave and satisfies $H_L(q,p)\leq -H_L(p, q)$. In most concrete examples, it is actually anti-symmetric.
If now $A$ is a maximal monotone operator, then $A^{-1}$ is also maximal monotone and therefore can be written as $A^{-1}=\bar \partial L$, where $L$ is a selfdual Lagrangian on $X^*\times X$ that can be constructed in the following way: First, let
\begin{equation}
N(p,x)=\sup\{\langle p,y\rangle +\langle q, x-y\rangle; \, (y,q)\in {\rm Graph}(A)\}
\end{equation}
in such a way that
 \begin{equation}
\hbox{$N^*(x,p) \geq N(p,x)\geq \langle x, p\rangle$ for every $(x,p)\in X\times X^*$.}
\end{equation}
Then consider the following Lagrangian on $X^*\times X$,
 \[
 L(x, p):=\inf \left\{\frac{1}{2}N(p_1, x_1)+\frac{1}{2}N^*(x_2, p_2)+\frac{1}{8}\|x_1-x_2\|^2+\frac{1}{8}\|p_1-p_2\|^2; \, (x, p)=\frac{1}{2}(x_1, p_1) + \frac{1}{2}(x_2, p_2)\right\}.
 \]
It was shown in \cite{Gh} that $L$ is a self-dual Lagrangian on $X^*\times X$. One can also show that the corresponding Hamiltonian $H_L$ on $X\times X$ is anti-symmetric, and that
\begin{equation}
\hbox{$p \in Ax$  \qquad \text{ if and only if } \qquad $(p,-p) \in \partial H_L (x,x)$.}
\end{equation}
Moreover, $Ax=\nabla_1 H(x,x)$ if $A$ is single-valued maximal monotone operator, which is Krause's result mentioned above.

Compared to Brenier's, our decomposition now looks like we have replaced the potential of a convex function in Brenier's theorem by a more general maximal monotone operator $A$ (or a self-dual Lagrangian $L$), while we have strengthened $S$ to make it a measure preserving involution. \\

\textbf{6. Connection to Monge's mass transport:} We shall see in the next section that the transformation $S$ appearing in the decomposition (\ref{rep.1}) of $u$ maximises the quantity
\[
\int_\Omega \langle u(x), S(x) \rangle \, dx,
\]
 among all measure preserving involutions on $\Omega$.  Equivalently, it does minimize
 \[
 \int_\Omega |u(x)-S(x)|^2 \, dx,
 \]
  which is the distance of $u$ to the set of measure preserving involutions on $\Omega$. The latter minimization can now be related to an optimal transport problem with a quadratic cost, between the measure $\mu$ on $\Omega \times u(\Omega)$ obtained by pushing Lebesgue measure on $\Omega$ by the map $x\to (x, u(x))$, and the measure $\nu$ on  $u(\Omega) \times \Omega$ obtained from $\mu$ by the transposition map $(x, y) \to (y,x)$. Indeed, any map $T$ pushing $\mu$ onto $\nu$ can be parametrized by an application $S:\Omega \to \Omega$ via the formula:
\[
T: (x, u(x)) \to (u(Sx), Sx),
\]
and the transport cost is then equal to
\[
\frac{1}{2}\int_\Omega [|u(Sx)-x|^2 +|u(x)-S(x)|^2] \, dx,
\]
which, in the case where $S$ is a measure preserving involution, coincides with
$
\int_\Omega |u(x)-S(x)|^2 \, dx.
$

Note now that if $T$ is an optimal transport mapping $\mu$ onto $\nu$, then the map $(y,x)\to (x,y) \to T(x,y)$ would be an optimal transport mapping $\nu$ onto $\mu$. It will then follow that if the optimal transport $T$ from $\mu$ onto $\nu$ was unique, then the $S$ corresponding to $T$ would necessarily be an involution. Now in terms of Brenier's theorem, the uniqueness would necessarily lead to a convex function $L$ on $\R^N\times \R^N$ such that $T=\nabla L$ and where $L$ is a selfdual Lagrangian, i.e., $L(x, p)=L^*(p,x)$. The anti-symmetric Hamiltonian $H$ would simply be the Legendre transform of $L$ with respect to the first variable.

Unfortunately, the measures $\mu$ and $\nu$ on the product space are too degenerate to fall under the framework where we have uniqueness in Brenier's theorem. Hence the need to establish the result directly and without resorting to Mass transport.

On the other hand, if one starts with a measure $\mu$ on the product space $\Omega \times \Omega$ that is absolutely continuous with respect to Lebesgue measure, then one can apply Brenier's theorem to find an optimal transport map $\nabla L$ that pushes $\mu$ onto its transpose $\tilde \mu$ via the involution $R(x,y)=(y,x)$, and where $L$ is a convex function on $\R^N\times \R^N$. Consider now the convex function $M(x,y)=L^*(y,x)=L^*\circ R (x,y)$ where $L^*$ is the Legendre transform of $L$ with respect to both variable, and note that $\nabla M=R\circ \nabla L^*\circ R$ also maps $\mu$ into $\tilde \mu$. By the uniqueness of the optimal map we deduce that $\nabla L=\nabla M= R\circ \nabla L^*\circ R$, which means that $L$ is a self-dual Lagrangian (up to a constant).

We are extremely grateful to Bernard Maurey for very insightful discussions, and to Ivar Ekeland and Philip Loewen for their helpful input.

\section{A variational formulation of the problem}
The  proof of the standard polar factorization by Brenier  was based on  tools from the  Monge-Kantorovich theory. Later W. Gangbo \cite{Gg} provided a more direct  proof by formulating the Brenier decomposition as the Euler-Lagrange equation corresponding to  a suitable variational problem.  Our approach is in line with Gangbo's method and involves various new results about skew-symmetric  functions, which may have their own interest.\\

To formulate our variational problem, we start by considering the set  ${\cal H}$ of all continuous anti-symmetric functions on $\bar \Omega$, that is
\begin{eqnarray*}
{\cal H}=\big \{H \in C(\bar \Omega \times \bar
\Omega);H(x,y)=-H(y,x) \}.
\end{eqnarray*}
For each $H \in {\cal H},$ define $L_{H}:\Omega \times \R^N\rightarrow \R$ by
\begin{equation}\label{LH}L_H(x,p)=\sup_{y \in \bar \Omega} \{ \langle y,p \rangle -H(y,x)\},\end{equation}
where $\langle .,.\rangle$ is the standard inner product in $\R^N.$ Set
$ {\cal L}=\big \{L_{H};
     H \in {\cal H} \big \}.$

  Here is our main theorem.
\begin{theorem}\label{main.2} Let $\Omega$ be an open bounded set in $\R^N$ such that $mea(\partial \Omega)=0$, and let $u \in L^{\infty}(\Omega, \R^N)$ be a non-degenerate vector field. Consider the following two variational problems:
\begin{eqnarray}\label{primal}
\hbox{$ P_{\infty}=\inf \big \{\int_{\Omega }L_{H}(x,u(x)) \, dx; H \in {\cal
 H} \big \}$}
   \end{eqnarray}
   and
   \begin{eqnarray}\label{dual}
\hbox{$ D_{\infty}=\sup \big \{\int_{\Omega }\langle u(x), S(x)\rangle \, dx; S:\bar
{\Omega}\rightarrow \bar {\Omega} \text{ is a measure preserving involution} \big \}.$}
   \end{eqnarray}
   Then the following assertions holds:
   \begin{enumerate}
  \item The variational problems (\ref{primal}) and (\ref{dual}) are dual
   in the sense that $P_{\infty}= D_{\infty}.$
   \item There exists a globally Lipschitz, anti-symmetric, convex-concave Hamiltonian $H$ on $\R^N \times \R^N$, such that the minimum in (\ref{primal}) is attained at  $L_{\bar H}$, where $\bar H$ is the restriction of $H$ to $\bar \Omega \times \bar \Omega$.

   \item  There exists a measure preserving involution $S$ such that
   the maximum in (\ref{dual}) is attained.
   \item  Moreover, for each $H$ satisfying (2), there exists $S$ satisfying (3) such that the following equation holds:
\begin{eqnarray}\label{rep.2}
u(x)= \nabla_1 H(S(x), x) \qquad \text{ a.e. } \quad x \in \Omega.
\end{eqnarray}

\end{enumerate}
\end{theorem}

 \subsection{Self-dual transformations}

 We first introduce the following notion.
 \begin{definition} \rm Let $\Omega$ be an open bounded subset of $\R^N$, and say that a measurable point transformation  $S:\bar
{\Omega}\rightarrow \bar {\Omega}$ is  {\it self-dual} if for every $H \in L^1(\Omega\times \Omega)$ such that $H(x,y)=-H(y,x)$ for almost all $(x,y)\in \Omega\times \Omega$, we have
\begin{equation}
\int_{\Omega}H(x, S(x)) \, dx=0.
\end{equation}
   \end{definition}
   In order to characterize these maps, we recall  the following basic notion.
    \begin{definition} \rm   A map $S:\bar
{\Omega}\rightarrow \bar {\Omega}$  is said to be a  {\it measure
preserving transformation} if for every $f \in L^1(\Omega)$, $f\circ S \in
L^1(\Omega)$ and
\begin{equation}
\int_{\Omega}f\circ S(x) \, dx=\int_{\Omega}f(x) \, dx.
\end{equation}
   \end{definition}

   By considering functions of the form $H(x,y)=f(x)-f(y)$ where $f \in L^1(\Omega)$, one can easily see that self-dual transformations are necessarily
measure preserving.

The converse is however not true. Indeed, the map $S(x_1, x_2)=(x_2, -x_1)$ on $R^2$  is clearly measure preserving. On the other hand, by considering the symplectic matrix $J(x_1, x_2)=(-x_2, x_1)$ on $R^2$, the function $H_J(x,y)=\langle Jx, y\rangle$ is clearly anti-symmetric. We therefore have that
\[
\int_{B}H_J(x, Sx)dx=\int_{B}\langle Jx, Sx\rangle dx=-\int_{B}\langle Jx, Jx\rangle dx= -\int_{B}\|x\|^2 dx \neq 0.
\]
which means that $S$  is not  a self-dual transformation.

More generally, recall that the linear and discrete analogue of measure preserving maps are the unitary matrices, i.e., those that satisfy $UU^*=U^*U=I$. If now $U$ is a self-dual transformation, then one can easily see that  $U^*=U$ and $U$ is therefore a symmetric involution. This turned out to be true for more general transformations.

\begin{proposition}\label{char} A measurable map $S$ is a self-dual point transformation on $\Omega$ if and only if it is both measure preserving and an involution, i.e., $S^2=I$ a.e., where $I$ is the identity map on $\Omega$.

\end{proposition}

\textbf{Proof.} Assuming that $S$ is measure preserving such that $S^2=I$ a.e, then  for every anti-symmetric
$H$ in $L^1(\Omega \times \Omega)$, we have
\begin{eqnarray*}
\int_{\Omega} H(x,S(x)) \, dx=\int_{\Omega} H(S(x), S^2(x)) \, dx=\int_{\Omega} H(S(x),x) \, dx
=-\int_{\Omega} H(x,S(x))\, dx,
\end{eqnarray*}
hence $\int_{\Omega} H(x,S(x)) \, dx=0$. That is $S$ is a self-dual transformation.\\

Conversely, Let $S$ be a self-dual transformation. We have seen above that it is necessarily measure preserving. Consider now the anti-symmetric functional  $H(x, y)=|S(x)-y|-|S(y)-x|$. We must have that
\[
 0=\int_{\Omega} H(x,S(x)) \, dx=\int |S^2(x)-x| dx,
 \]
which clearly yields that $S$ is an involution almost everywhere.

An immediate application of Proposition \ref{char} is that
\[
P_\infty \geq D_\infty.
\]
Indeed, for any $H\in {\cal H}$, and any point transformation $S$ on $\Omega$, we have
\[
L_H(x, u(x))+H(S(x), x) \geq \langle u(x), S(x)\rangle.
\]
 If now $S$ is a measure preserving involution, we have that $\int_{\Omega} H(x,S(x)) \, dx=0$, which means that
 \[
 \int_\Omega L_H(x, u(x)) \geq \int_\Omega \langle u(x), S(x)\rangle \, dx.
 \]

\subsection{Regularization of skew-symmetric functions}
 Let $\Omega$ be a bounded domain contained in a ball $B_{R}$ centered at the origin with radius $R>0$ in $\R^N$, and consider an arbitrary anti-symmetric function $H \in C(\bar \Omega \times \bar
\Omega).$  Our plan is to show that one can assign to $H,$ a skew-symmetric convex-concave Lipschitz function $H_{reg}$ such that $L_{H_{reg}} \leq L_{H}$ on $\bar \Omega \times B_R.$

 Note that for  an elementary anti-symmetric function of the form $H(x,y)=f(x)-f(y)$, where $f$ is a  continuous function on $\bar \Omega,$ one can easily show that $H_{reg}(x,y)=f^{**}(x)-f^{**}(y)$ where $f^{**}$ is a double conjugate of $f$ defined as follows:
 \[
\hbox{$ f^*(p)=\sup\limits_{y \in \bar \Omega} \{ \langle x,p \rangle -f(y)\}$ and
$f^{**}(x)=\sup\limits_{p \in \bar B_R} \{ \langle x,p \rangle -f^*(p)\}$.}
\]
Note that these are not the usual Legendre transforms since the suprema are not taken over the whole space.
 The analogous regularization process for a general anti-symmetric function is more technical.

 We start by considering the class
\begin{eqnarray}
{\cal H}_-=\big \{H \in C(\bar \Omega \times \bar
\Omega);H(x,y)\leq -H(y,x)  \, \text{ for all } x, y \in \Omega\}.
\end{eqnarray}
For each $H \in {\cal H}_-$ define $L_H$ as in (\ref{LH}) and set
\[
{\cal L}_-=\{L_H;H \in {\cal H}_- \}.
\]
 For $L_{H} \in {\cal L}_-$,  we define its (restricted) Fenchel dual $L^*_{H}:\R^N \times \R^N \to \R$ by
\begin{eqnarray*}
L_{H}^*(q,y)=\sup_{x\in \bar \Omega, p \in B_R}\{\langle y,p\rangle+
\langle q,x \rangle- L_{H}(x,p)\}.
\end{eqnarray*}
and similarly $L^{**}_{H}:\R^N \times \R^N \to \R$ by
\begin{eqnarray*}
L_{H}^{**}(y,q)=\sup_{x\in \bar \Omega, p \in B_R}\{\langle y,p\rangle+
\langle q,x \rangle- L^*_{H}(p,x)\}.
\end{eqnarray*}
For each convex function  $L:\R^N \times \R^N \to \R$, we shall define its $B_R$-Hamiltonian by
\[H_{L}(x,y)= \sup_{ p \in B_R}\{\langle x,p\rangle- L(y,p)\}\]

Finally for each $H \in {\cal H},$ we define the convex-concave regularized $H_{reg}$ of $H$ by
\begin{eqnarray*}
H_{reg}(x,y)=\frac{H_{L_H^{**}}(x,y)-H_{L_H^{**}}(y,x)}{2}, \qquad \qquad \text{ for all } x,y \in \R^N.
\end{eqnarray*}
We list some of the properties of $H_{reg}$ and  $L_{ H_{reg}}.$

\begin{proposition} \label{prop1} Let   $H \in {\cal H}.$ The following statements hold:
\begin{enumerate}
\item   $ H_{reg}$ is a skew-adjoint Hamiltonian on $\R^N\times \R^N$ whose restriction to $\bar \Omega\times \bar \Omega$ belong to$ {\cal
H}$.
\item   $L_{ H_{reg}}$ is convex and continuous in both variables and  $L_{ H_{reg}} \leq L_{ H}$ on $\bar \Omega \times B_R$.
\item   $|L_{ H_{reg}}(x,p)| \leq R\|x\|+R\|p\|+5R^2$ and $| H_{reg}(x,y)| \leq R\|x\|+R\|y\|+4R^2$ for all $x,y,p \in \R^N.$
\item    $L_{ H_{reg}}$ and $H_{reg}$ are both Lipschitz continuous with Lipschitz constants less than
$4NR.$
\end{enumerate}
\end{proposition}
To prove the above Proposition, we shall need the following two lemmas.
\begin{lemma}\label{ine1} If $H \in {\cal H}_-$, then hen $H_{L^{**}_{H}} \in {\cal H}_-.$
\end{lemma}

\textbf{Proof.} For $x,y \in \R^N$ we have
\begin{eqnarray*}
H_{L_H^{**}} (x, y)=\sup_{ p_1 \in B_R}\{ \langle p_1,x \rangle- L_H^{**}(y, p_1)\}.
\end{eqnarray*}
It follows from the definition of $L_H^{**}$ that
\begin{eqnarray}\label{l88}
L_H^{**}(y, p_1)&=&\sup_{z_1 \in \bar
\Omega, p_2 \in B_R}\Big \{\langle z_1,p_1\rangle+ \langle
p_2,y \rangle- L^*_H(p_2,z_1)\Big\} \nonumber \\
&=&\sup_{z_1\in \bar \Omega, p_2 \in B_R}\Big \{
\langle z_1,p_1\rangle+\langle
p_2,y \rangle-\sup_{z_2\in \bar \Omega, p_3 \in B_R}\{\langle z_2,p_2\rangle+ \langle
p_3,z_1 \rangle- L_H(z_2, p_3)\Big\} \nonumber\\
&=&\sup_{z_1\in \bar \Omega, p_2 \in B_R}\inf_{z_2\in \bar \Omega, p_3 \in B_R}\Big \{
\langle z_1,p_1\rangle+ \langle
p_2,y \rangle-\langle z_2,p_2\rangle- \langle
p_3,z_1 \rangle  \nonumber \\ &&+\sup_{z_3\in \bar \Omega}\{ \langle z_3,p_3\rangle-H(z_3, z_2)\} \Big\}  \nonumber\\
&=&\sup_{z_1\in \bar \Omega, p_2 \in B_R}\inf_{z_2\in \bar \Omega, p_3 \in B_R}\sup_{z_3\in \bar \Omega}\Big \{
\langle z_1,p_1\rangle+ \langle
p_2,y \rangle-\langle z_2,p_2\rangle- \langle
p_3,z_1 \rangle+ \langle z_3,p_3\rangle -H(z_3, z_2) \Big\}. \nonumber\\
\end{eqnarray}

Therefore,
\begin{eqnarray}\label{LH1}
H_{L_H^{**}} (x, y)&=&\sup_{ p_1 \in B_R}\inf_{z_1\in \bar \Omega, p_2 \in B_R}\sup_{z_2\in \bar \Omega, p_3 \in B_R}\inf_{z_3\in \bar \Omega}\Big \{
\langle p_1,x \rangle-\langle z_1,p_1\rangle- \langle
p_2,y \rangle+\langle z_2,p_2\rangle+ \langle
p_3,z_1 \rangle- \langle z_3,p_3\rangle \nonumber \\&&+H(z_3, z_2) \Big\}.
\end{eqnarray}
Taking into account that  $H(z_3, z_2) \leq -H(z_2, z_3)$ together with the fact that
\[
\sup_{z_2 \in \bar \Omega}\inf _{z_3 \in \bar \Omega} \{...\} \leq \inf _{z_3 \in \bar \Omega} \sup_{z_2 \in \bar \Omega }\{...\}
\]
 yield
\begin{eqnarray*}
H_{L_H^{**}} (x, y)&\leq&\sup_{ p_1 \in B_R}\inf_{z_1\in \bar \Omega, p_2 \in B_R}\sup_{ p_3 \in B_R} \inf_{z_3\in \bar \Omega}\sup_{z_2\in \bar \Omega}\Big \{
\langle p_1,x \rangle-\langle z_1,p_1\rangle- \langle
p_2,y \rangle+\langle z_2,p_2\rangle+ \langle
p_3,z_1 \rangle- \langle z_3,p_3\rangle\\&&-H(z_2, z_3) \Big\}\\
&=&\sup_{ p_1 \in B_R}\inf_{z_1\in \bar \Omega, p_2 \in B_R}\sup_{ p_3 \in B_R} \inf_{z_3\in \bar \Omega}\Big \{
\langle p_1,x \rangle-\langle z_1,p_1\rangle- \langle
p_2,y \rangle+ \langle
p_3,z_1 \rangle- \langle z_3,p_3\rangle\\&&\sup_{z_2\in \bar \Omega}\{ \langle z_2,p_2\rangle-H(z_2, z_3)\} \Big\}.\\
&=&\sup_{ p_1 \in B_R}\inf_{z_1\in \bar \Omega, p_2 \in B_R}\sup_{ p_3 \in B_R} \inf_{z_3\in \bar \Omega}\Big \{
\langle p_1,x \rangle-\langle z_1,p_1\rangle- \langle
p_2,y \rangle+ \langle
p_3,z_1 \rangle- \langle z_3,p_3\rangle+L_H(z_3,p_2)\Big\}
\end{eqnarray*}
Note also that $\inf_{z_3\in \bar \Omega}\{...\} \leq \{...\}_{|z_3=z_1},$ from which we obtain
\begin{eqnarray*}
H_{L_H^{**}} (x, y)&\leq& \sup_{ p_1 \in B_R}\inf_{ p_2 \in B_R}\inf_{z_1\in \bar \Omega} \sup_{ p_3 \in B_R} \Big \{
\langle p_1,x \rangle-\langle z_1,p_1\rangle- \langle
p_2,y \rangle+ \langle
p_3,z_1 \rangle- \langle z_1,p_3\rangle+L_H(z_1,p_2)\Big\}\\
&=& \sup_{ p_1 \in B_R}\inf_{ p_2 \in B_R}\inf_{z_1\in \bar \Omega}\sup_{ p_3 \in B_R} \Big \{
\langle p_1,x \rangle-\langle z_1,p_1\rangle- \langle
p_2,y \rangle+L_H(z_1,p_2)\Big\}.
\end{eqnarray*}
We  can now drop the term $\sup_{ p_3 \in B_R}$ since there is no $p_3$ in the expression in the bracket. Thus,
 \begin{eqnarray*}
H_{L_H^{**}} (x, y)&\leq& \sup_{ p_1 \in B_R}\inf_{ p_2 \in B_R}\inf_{z_1\in \bar \Omega}\Big \{
\langle p_1,x \rangle-\langle z_1,p_1\rangle- \langle
p_2,y \rangle+L_H(z_1,p_2)\Big\}\\
&=& \sup_{ p_1 \in B_R}\Big \{
\langle p_1,x \rangle +\inf_{ p_2 \in B_R} \inf_{z_1\in \bar \Omega}\{-\langle z_1,p_1\rangle- \langle
p_2,y \rangle+L_H(z_1,p_2)\}\Big\}\\
&=& \sup_{ p_1 \in B_R}\Big \{
\langle p_1,x \rangle -L_H^*(p_1, y)\Big\}.
\end{eqnarray*}
It follows from the definition of $L^{**}_{H}$ that $L_H^*(p_1,y)+L^{**}_{H}(z,p)\geq \langle p_1,z \rangle+\langle y,p\rangle$ for all $z, p \in \bar \Omega \times  B_R.$ Therefore
\begin{eqnarray*}
H_{L_H^{**}} (x, y)&\leq& \sup_{ p_1 \in B_R}\Big \{
\langle p_1,x \rangle -L_H^*(p_1, y)\Big\}\\
&\leq& \sup_{ p_1 \in B_R}\Big \{
\langle p_1,x \rangle -\sup_{ p \in B_R, z\in \bar \Omega} \{\langle p_1,z \rangle+\langle y,p\rangle-L_H^{**}(z,p)\}\Big\}\\
&=& \sup_{ p_1 \in B_R}\inf_{ p \in B_R, z\in \bar \Omega}\Big \{
\langle p_1,x \rangle -\langle p_1,z \rangle-\langle y,p\rangle+L_H^{**}(z,p)\Big\}\\
&\leq& \inf_{ p \in B_R, z\in \bar \Omega}\sup_{ p_1 \in B_R}\Big \{
\langle p_1,x \rangle -\langle p_1,z \rangle-\langle y,p\rangle+L_H^{**}(z,p)\Big\}\\
&\leq& \inf_{ p \in B_R}\sup_{ p_1 \in B_R}\Big \{
\langle p_1,x \rangle -\langle p_1,x \rangle-\langle y,p\rangle+L_H^{**}(x,p)\Big\}\\
&=&\inf_{ p \in B_R}\Big \{
-\langle y,p\rangle+L_H^{**}(x,p)\}\Big\}\\
&=&-H_{L_H^{**}} (y, x).
\end{eqnarray*}
This proves that $H_{L_H^{**}} (x, y) \leq -H_{L_H^{**}} (y, x)$ for all $x, y \in \bar \Omega.$
\hfill $\square$

\begin{lemma}\label{ine2} If $H \in {\cal H}_-$, then  $L^*_H(p,x) \leq L_H(x,p)$  and $L^{**}_H(x,p) \leq L_H(x,p)$ for all $x\in \bar \Omega, p \in B_R.$
\end{lemma}
\textbf{Proof.} For every $(x,p) \in \bar \Omega \times B_R,$ we have
\begin{eqnarray}\label{LH*}
L_H^*(p,x)&=&\sup_{y \in \bar \Omega, q \in B_R} \{ \langle y,p
\rangle+\langle x, q \rangle -L_{H}(y,q)\}\nonumber \\
&=&\sup_{y \in \bar
\Omega, q \in B_R} \{ \langle y,p
\rangle+\langle x, q \rangle -\sup_{z \in \bar \Omega}\{\langle q,z\rangle-H(z,y)\}\} \nonumber
\\&=&\sup_{y \in \bar \Omega,
q \in B_R} \inf_{z \in \bar \Omega}\{ \langle y,p
\rangle+\langle x, q \rangle -\langle q,z\rangle+H(z,y)\}.
\end{eqnarray}
It then follows that
\begin{eqnarray*}
L_H^*(p,x)&\leq &\sup_{y \in \bar \Omega,
q \in B_R} \inf_{z \in \bar \Omega}\{ \langle y,p
\rangle+\langle x, q \rangle -\langle q,z\rangle-H(y,z)\}\\
\\&\leq&\inf_{z \in \bar \Omega} \sup_{y \in \bar \Omega,
q \in B_R} \{ \langle y,p
\rangle+\langle x, q \rangle -\langle q,z\rangle-H(y,z)\}\\
&=&\inf_{z \in \bar \Omega} \sup_{y \in \bar \Omega} \{ \langle y,p
\rangle+R\|x-z\|-H(y,z)\}\\
&=&\inf_{z \in \bar \Omega}  \{ R\|x-z\|+L_H(z,p)\}\\
&\leq& L_H(x,p).
\end{eqnarray*}
Thus, $L_H^*(p,x) \leq L_H(x,p).$  We now prove $L^{**}_H\leq L_H$ on $\bar \Omega \times B_R.$
It follows from equality (\ref{l88}) that for every $(y,p_1) \in \bar \Omega \times B_R,$
\begin{eqnarray*}
L_H^{**}(y, p_1)
=\sup_{z_1\in \bar \Omega, p_2 \in B_R}\inf_{z_2\in \bar \Omega, p_3 \in B_R}\sup_{z_3\in \bar \Omega}\Big \{
\langle z_1,p_1\rangle+ \langle
p_2,y \rangle-\langle z_2,p_2\rangle- \langle
p_3,z_1 \rangle+ \langle z_3,p_3\rangle -H(z_3, z_2) \Big\}.
\end{eqnarray*}
Thus,
\begin{eqnarray*}
L_H^{**}(y, p_1)
&\leq & \inf_{z_2\in \bar \Omega, p_3 \in B_R}\sup_{z_1\in \bar \Omega, p_2 \in B_R}\sup_{z_3\in \bar \Omega}\Big \{
\langle z_1,p_1\rangle+ \langle
p_2,y \rangle-\langle z_2,p_2\rangle- \langle
p_3,z_1 \rangle+ \langle z_3,p_3\rangle -H(z_3, z_2) \Big\}
\end{eqnarray*}
Therefore by considering $z_2=y$ and $p_3=p_1$ we get
\begin{eqnarray*}
L_H^{**}(y, p_1)
\leq  \sup_{z_3\in \bar \Omega}\{
  \langle z_3,p_1\rangle -H(z_3, y)\}
=L_H(y, p_1).
\end{eqnarray*}
This completes the proof.
\hfill $\square$ \\

We now recall the following result from \cite{Gg}.
\begin{proposition}\label{Gp} Let $D$ be an open set in $\R^d$ such that $\bar D \subset \tilde B_R$ where $\tilde B_R$ is ball with radious $R$ centered at the origin in $\R^d.$ Let $f:\R^d \to \R$ and define $\tilde f : \R^d \to \R$  by
\[\tilde f(y)= \sup_{z \in D}\{\langle y,z\rangle-f(z)\}.\]
If $f \in L^{\infty}(D),$ then $\tilde f$ is a convex, Lipschitz function and
\[|\tilde f (y_1)-\tilde f(y_2)| \leq dR \|y_1-y_2\|,\]
for all $y_1, y_2 \in \R^d.$
\end{proposition}
\begin{lemma} \label{prop2} If $H \in {\cal H}_-$, then the following statements hold:
\begin{enumerate}
\item  $|L^{**}_H(x,p)| \leq R\|x\|+R\|p\|+3R^2 $ and $|H_{L^{**}_H}(x,y)| \leq R\|x\|+R\|y\|+4R^2 $ for all $x,y, p \in \R^N.$
\item $L^{**}_{H}$ and  $H_{L^{**}_H}$ are Lipschitz with $Lip(H_{L^{**}_H}), Lip(L^{**}_H)\leq 4NR.$
\end{enumerate}
\end{lemma}
\textbf{Proof.} It follows from the definition of $L_H$ that $L_{H}(y,q) \geq  \langle y,q\rangle$ on $\bar \Omega \times B_R$. This together with $\bar \Omega \subset B_R$ imply that
\begin{eqnarray*}
L_{H}^*(p,x)&=&\sup_{y\in \bar \Omega, q \in B_R}\{\langle x,q\rangle+
\langle p,y \rangle- L_{H}(y,q)\}\\
&\leq & \sup_{y\in \bar \Omega, q \in B_R}\{\langle x,q\rangle+
\langle p,y \rangle- \langle y,q\rangle\}\\
&\leq & R\|x\|+R\|p\|+R^2
\end{eqnarray*}
for all $x, p \in \R^N.$
With a similar argument we obtain $L^{**}_H(x,p) \leq R\|x\|+R\|p\|+R^2$. We also have
\begin{eqnarray*}
L_{H}^{**}(x,p)&=&\sup_{y\in \bar \Omega, q \in B_R}\{\langle x,q\rangle+
\langle p,y \rangle- L^*_{H}(y,q)\}\\
&\geq & \langle x,q\rangle+
\langle p,y \rangle- L^*_{H}(y,q), \qquad \qquad (\text{ for all } (y,q) \in \bar \Omega \times  B_R)\\
&\geq &- R\|x\|-R\|p\|-R\|y\|-R\|q\|-R^2\\ &\geq& - R\|x\|-R\|p\|-3R^2.
\end{eqnarray*}
Therefore $ |L^{**}_H(x,p)| \leq R\|x\|+R\|p\|+3R^2.$ The estimate for $H_{L^{**}_H}$  can be easily deduced from its definition together with the  estimate on $L^{**}_H.$ This completes the proof of part (1).\\

For (2) set $D=  \Omega \times B_R$, then $D \subset \tilde B_{2R}$ where  $\tilde B_{2R}$ is a ball with radius $2R$ in  $\R^{2N}.$ Now assuming $f=L^*_H$ in Proposition \ref{Gp}, we have that $\tilde f =L^{**}_H$. Therefore $L^{**}_H$ is Lipschitz in $\R^{2N}$ with $Lip(L^{**}_H) \leq 4NR.$ To prove that $H_{L^{**}_H}$ is Lipschitz we first fix $y \in \R^{N}$ and we define $f_y: \R^{N} \to \R$ by $f_y(p)=L^{**}_H(y,p)$. Assuming $D=B_R \subset \R^{N}$ in Proposition \ref{Gp}, we obtain that the map $x \to \tilde f_y(x)=H_{L^{**}_H}(x,y)$ is Lipschitz and
\begin{equation}\label{lip}
|H_{L^{**}_H}(x_1,y)-H_{L^{**}_H}(x_2,y)| \leq NR \|x_1-x_2\|
\end{equation}
for all $x_1, x_2 \in \R^{N}.$ Noticing that the Lipschitz constant $NR$ is independent of $y,$ the above inequality holds for all $x_1, x_2, y \in \R^{N}.$ To prove $H_{L^{**}_H}(x,y)$ is Lipschitz with respect to the second variable,  let  $r>0$ and $y_1, y_2 \in \R^{N}.$ Let $p_1$ and  $p_2$ be such that
\[\langle x,p_2\rangle- L^{**}_{H}(y_1,p_2) \leq H_{L^{**}_H}(x,y_1) \leq \langle x,p_1\rangle- L^{**}_{H}(y_1,p_1)+ r,\]
and
\[\langle x,p_1\rangle- L^{**}_{H}(y_2,p_1) \leq H_{L^{**}_H}(x,y_2) \leq \langle x,p_2\rangle- L^{**}_{H}(y_2,p_2)+ r,\]
It follows that
\[L^{**}_{H}(y_2,p_2)-L^{**}_{H}(y_1,p_2)- r\leq   H_{L^{**}_H}(x,y_1)-H_{L^{**}_H}(x,y_2) \leq L^{**}_{H}(y_2,p_1) -L^{**}_{H}(y_1,p_1) +r\]
$L^{**}_{H}$ is Lipschitz, thus,
\[-4NR\|y_1-y_2\|- r\leq   H_{L^{**}_H}(x,y_1)-H_{L^{**}_H}(x,y_2) \leq 4NR\|y_1-y_2\| +r\]
  Since $r>0$ is arbitrary we obtain \[-4NR\|y_1-y_2\|\leq   H_{L^{**}_H}(x,y_1)-H_{L^{**}_H}(x,y_2) \leq 4NR\|y_1-y_2\|.\]
This together with ( \ref{lip}) prove that $H_{L^{**}_H}$ is Lipschitz and $Lip(H_{L^{**}_H})\leq 4NR.$
 \hfill $\square$\\

\textbf{Proof of Proposition \ref{prop1}.}
 It is  easily  seen  that $H_{reg}(x,y)=-H_{reg}(y,x)$ for all $x,y \in \R^N$.\\
Now note that by definition \[H_{L_H^{**}}(x,y)= \sup_{ p \in B_R}\{\langle x,p\rangle- L_H^{**}(p,y)\},
\]
and therefore for all $y \in \R^N,$ the function $x \to H_{L_H^{**}}(x,y)$ is convex. We shall show that  for all $x\in \R^N$ the function
$ y \to H_{L_H^{**}}(x,y)$ is concave. In fact we  need to show that \[y \to -H_{L_H^{**}}(x,y)=\inf_{ p \in B_R}\{L_H^{**}(p,y)-\langle x,p\rangle\}\]
is convex. For that  consider $ \lambda  \in (0,1)$ and elements $y_1, y_2 \in \R^N.$ For every  $a > -H_{L_H^{**}}(x,y_1)$ and $b > -H_{L_H^{**}}(x,y_2),$ find $p_1, p_2 \in \R^N$ such that
\[-H_{L_H^{**}}(x,y_1) \leq L_H^{**}(p_1,y_1)-\langle x,p_1\rangle \leq a \quad \text{ and }  -H_{L_H^{**}}(x,y_2) \leq L_H^{**}(p_2,y_2)-\langle x,p_2\rangle \leq b.\]
Now use the convexity of the ball $B_R$ and the convexity of the function $L_H^{**}$ in both variables  to write
\begin{eqnarray*}
-H_{L_H^{**}}(x,\lambda y_1+(1-\lambda) y_2)&=& \inf_{ p \in B_R}\{L_H^{**}(p,\lambda y_1+(1-\lambda) y_2 )-\langle x,p\rangle\}\\ & \leq & L_H^{**}(\lambda p_1+(1-\lambda)p_2,\lambda y_1+(1-\lambda) y_2 )-\langle x,\lambda p_1+(1-\lambda)p_2\rangle\}\\
& \leq & \lambda \big (L_H^{**}( p_1,y_1)-\langle x,p_1\rangle \big ) +(1-\lambda) \big (L_H^{**}( p_2,y_2)-\langle x,p_2\rangle \big )\}\\
& \leq & \lambda a +(1-\lambda)b.
\end{eqnarray*}
This establishes the of concavity of $y \to H_{L_H^{**}}(x,y)$. It then follows that $H_{reg}(.,y)$ is convex and $H_{reg}(x,.)$ is concave on $\R^N$,  which completes the proof of part (1).\\
We now prove part (2). Let $(x,p) \in \bar \Omega \times B_R.$ We have
\begin{eqnarray*}
L_{H_{reg}}(x,p)&=&\sup_{y\in \bar \Omega}\{\langle y,p\rangle- H_{reg}(y,x)\}\\
&=&\sup_{y\in \bar \Omega}\{\langle y,p\rangle- \frac{H_{L_H^{**}}(y,x)-H_{L_H^{**}}(x,y)}{2}\}\\
& \leq & \sup_{y\in \bar \Omega}\{\langle y,p\rangle- H_{L_H^{**}}(y,x)\} \qquad \quad \big (\text{by Lemma \ref{ine1}  } H_{L_H^{**}}(y,x) \leq -H_{L_H^{**}}(x,y) \big )\\
&=& \sup_{y\in \bar \Omega}\big \{\langle y,p\rangle- \sup_{q\in \bar B_R}\{\langle y,q\rangle- L_H^{**}(x,q)\}\big \}\\
&=& \sup_{y\in \bar \Omega}\inf_{q\in \bar B_R}\big \{\langle y,p\rangle- \langle y,q\rangle+ L_H^{**}(x,q)\big \}\\
&\leq & \inf_{q\in \bar B_R}\sup_{y\in \bar \Omega}\big \{\langle y,p\rangle- \langle y,q\rangle+ L_H^{**}(x,q)\big \}\\
&\leq & \sup_{y\in \bar \Omega}\big \{\langle y,p\rangle- \langle y,p\rangle+ L_H^{**}(x,p)\big \}\\
&=& L_H^{**}(x,p).
\end{eqnarray*}
Thus, $L_{H_{reg}}(x,p) \leq L_H^{**}(x,p)$. It also follows from Lemma \ref{ine2} that $  L_H^{**}(x,p) \leq  L_H(x,p)$ from which we obtain $L_{ H_{reg}} \leq L_{ H}$ on $\bar \Omega \times B_R$.\\

By the definition
$H_{reg}(x,y)=\frac{H_{L_H^{**}}(x,y)-H_{L_H^{**}}(y,x)}{2}.$ Therefore, by  Lemma \ref{prop2}, $H_{reg}$ is Lipschitz and $Lip(H_{reg}))=Lip(H_{L_H^{**}}) \leq 4NR$ and also \[|H_{reg}(x,y)| \leq \frac{|H_{L_H^{**}}(x,y)|+|H_{L_H^{**}}(y,x)|}{2} \leq R\|x\|+R\|y\|+4R^2.\] The corresponding results for $L_{H_{reg}}$ follow
 by the same arguments as in parts (2) of Lemma \ref{prop2} and the bound for $H_{reg}$. \hfill $\square$\\

\subsection{Proof of Theorem \ref{main.2}}

We first show that the minimization problem (\ref{primal}) has a solution. Let $B_R$ be a ball such that $\bar \Omega, u(\bar \Omega) \subset B_R.$
 Let $\{H^n\}$ be a sequence in ${\cal H}$ such that $L_{H^n}$ is a minimizing sequence for $P_{\infty}.$ It follows from Proposition \ref{prop1} that $L_{H_{reg}^n} \leq L_{H^n}$ on $\bar \Omega \times B_R$ and therefore $L_{H_{reg}^n}$  is still a minimizing for $P_{\infty}.$ It also follows from Proposition \ref{prop1} that  $L_{H_{reg}^n}$ and $H_{reg}^n$
are uniformly  Lipschitz with  $Lip(H^n_{reg}), Lip(L_{H^n_{reg}}) \leq 4NR$ and also \[|H^n_{reg}(x,y)| \leq R\|x\|+R\|y\|+4R^2 \quad \text{ and } \quad |L_{H^n_{reg}}(x,p)| \leq R\|x\|+R\|p\|+5R^2 \]
  for all $x,y, p \in \R^N.$ By Arzela-Ascoli's theorem, there exists two Lipschitz functions $H, L: \R^N \times \R^N \to \R$ such that $H^n_{reg}$ converges to $H$ and $L_{H^n_{reg}}$ converges to $ L$ uniformly on every compact set of $\R^N \times \R^N.$  This implies that $H \in {\cal H}.$ Note that
\[L_{H^n_{reg}}(x,p) +H^n_{reg}(y,x) \geq  \langle y,p\rangle, \]
for all $x,p \in \R^N $ and $y \in \bar \Omega,$ from which we have
 \[L(x,p) \geq  \langle y,p\rangle -H(y,x) , \]
for all $x,p \in \R^N $ and $y \in \bar \Omega.$
It implies that $L_H \leq L$. Let $H_{reg}$ be the regularization of $H$ defined in the previous section. Set $H_{\infty}= H_{reg}$ and $L_{\infty}=L_{H_{\infty}}$. It  follows from Proposition \ref{prop1} that $L_{H_{\infty }} \leq L_H$ on $\bar \Omega \times B_R,$ from which we have
\[P_{\infty}=\int_{\Omega} L_H(x, u(x))\, dx= \int_{\Omega} L_{\infty}(x, u(x))\, dx.\]
\hfill $\square$

For the rest of the proof, we shall need the following two technical lemmas.

\begin{lemma} \label{wei} Let $H_\infty$ be the Hamiltonian obtained above.  For each $x \in \bar \Omega$, define $f_x: \R^N \to \R$ by $f_x(y)=H_{\infty}(y,x).$ We also define $\tilde f_x: \R^N \to \R \cup\{+\infty\}$ by
 $
\tilde f_x (y)=f_x(y)$ whenever $y \in \bar \Omega$ and $+\infty$ otherwise. Let $(\tilde f_x)^{*}$ be the standard Fenchel dual of $\tilde f_x$ on $\R^N$, in such a way that  $(\tilde f_x)^{***}=(\tilde f_x)^{*}$ on $\R^N.$ Then we have,
\begin{enumerate}
\item $f_x = (\tilde f_x)^{**}= \tilde f_x$ on $\bar \Omega$, and
\item
$L_{\infty}(x,p)=
\sup\limits_{z \in \bar \Omega} \{\langle z,p \rangle-(\tilde f_x)^{**}(z)\}=\sup\limits_{z \in \R^N} \{\langle z,p \rangle-(\tilde f_x)^{**}(z)\}.$
\end{enumerate}
\end{lemma}
\textbf{Proof.}
(1) Since $(\tilde f_x)^{**}$ is the largest convex function below $\tilde f_x$ we have  and $f_x \leq (\tilde f_x)^{**}\leq \tilde f_x,$ from which we obtain  $f_x = (\tilde f_x)^{**}= \tilde f_x$ on $\bar \Omega.$

For (2), we first deduce from (1) that
\begin{eqnarray*}
(\tilde f_x)^{*}(y)= (\tilde f_x)^{***}(y)&=& \sup_{z \in \R^N} \{\langle z,y \rangle-(\tilde f_x)^{**}(z)\}\\
&\geq& \sup_{z \in B_R} \{\langle z,y \rangle-(\tilde f_x)^{**}(z)\}\\
&\geq& \sup_{z \in \Omega} \{\langle z,y \rangle-(\tilde f_x)^{**}(z)\}\\
&=& \sup_{z \in \Omega} \{\langle z,y \rangle-f_x(z)\}\\
&=& \sup_{z \in \Omega} \{\langle z,y \rangle- \tilde f_x(z)\}\\
&=& (\tilde f_x)^{*}(y),
\end{eqnarray*}
from which we have the desired result.
 \hfill $\square$

\begin{lemma} \label{limit} Let $H \in {\cal H}.$ For each $r \in \R$ and $\lambda >0$ define
\begin{eqnarray*}
L_{r,\lambda}(x,p)&:=& \sup_{z \in \bar \Omega}\{\langle z,p
\rangle-(\tilde f_x)^{**}(z) -\lambda\frac{ \|z\|^2}{2}+\lambda\frac{\|x\|^2}{2} -r
H(z,x)\},\\
L_{\lambda}(x,p)&:=& \sup_{z \in \R^N}\{\langle z,p \rangle-(\tilde f_x)^{**}(z)- \lambda\frac{\|z\|^2}{2}+\lambda\frac{\|x\|^2}{2}\},\\
L_{r}(x,p)&:=& \sup_{z \in \bar \Omega}\{\langle z,p
\rangle-H_{\infty}(z,x)  -r
H(z,x)\}.
\end{eqnarray*}
Then the following assertions hold:
\begin{enumerate}
\item For every $(x, p)\in \R^N\times \R^N$, we have $\lim\limits_{\lambda \to 0^+} L_{\lambda}(x,p)=L_{\infty}(x,p)$  and $\lim\limits_{\lambda \to 0^+} L_{r,\lambda}(x,p)=L_{r}(x,p).$
\item  For all $x\in \R^N$, the function $p \to L_{\lambda}(x,p)$ is differentiable.
\item For every $(x, p)\in \R^N\times \R^N$, we also have
$\lim\limits_{r\rightarrow 0} \frac {L_{r, \lambda}(x,p)-L_{\lambda}(x,p)}{r}=H(\nabla_2
L_{\lambda}(x,p),x).$
\end{enumerate}
\end{lemma}
\textbf{Proof.} Yosida's regularization of convex functions and (1) of Lemma \ref{wei}  yield that \[\lim_{\lambda \to 0^+} L_{r,\lambda}(x,p)=\sup_{z \in \bar \Omega}\{\langle z,p
\rangle-(\tilde f_x)^{**}(z) -r
H(z,x)\}=\sup_{z \in \bar \Omega}\{\langle z,p
\rangle-H_{\infty}(z,x) -r
H(z,x)\}=L_{r}(x,p).\]
We also have
\[\lim_{\lambda \to 0}L_{\lambda}(x,p)= \sup_{z \in \R^N}\{\langle z,p \rangle-(\tilde f_x)^{**}(z)\},\]
which, together with (2) of Lemma \ref{wei}, yield $\lim_{\lambda \to 0}L_{\lambda}(x,p)=L_{\infty}(x,p).$

 (2) follows from the fact that the Yosida regularization of  convex functions are differentiable.

 (3) We let  $z_{r, \lambda} \in \bar \Omega$ and $ z'_{r, \lambda} \in \R^N$ be such that
\begin{eqnarray*}
L_{r,\lambda}(x,p)&\leq& \langle z_{r, \lambda},p \rangle-(\tilde f_x)^{**}(z_{r, \lambda})-\lambda\frac{ \|z_{r, \lambda}\|^2}{2}+\lambda\frac{\|x\|^2}{2}-r
H(z_{r, \lambda},x)+r^2,\\
L_{\lambda}(x,p)&\leq& \langle z'_{r, \lambda},p \rangle-(\tilde f_x)^{**}(z'_{r, \lambda})-\lambda\frac{ \|z'_{r, \lambda}\|^2}{2}+\lambda\frac{\|x\|^2}{2}+r^2.
\end{eqnarray*}
Therefore,
\begin{eqnarray}\label{ine}
-H(z'_{r, \lambda},x)-r\leq \frac {L_{r,\lambda}(x,p)-L_{\lambda}(x,p)}{r}\leq -H(z_{r, \lambda}, x)+r.
\end{eqnarray}

 Note that by the definition of  $L_\lambda,$ we have  $\sup_{r \in [-1,1]}\|z'_{r, \lambda}\| < \infty.$ Suppose now that, up to a subsequence, $z_{r, \lambda}\rightarrow z_{\lambda} \in \bar \Omega$ and $z'_{r, \lambda}\rightarrow z_{\lambda}'$ as $r \to 0.$  This together with the definition of
$L_{r,\lambda}$ and $L_{\lambda}$ imply that
\begin{eqnarray}\label{equality}
L_{\lambda}(x,p)= \langle z_{\lambda},p \rangle-(\tilde f_x)^{**}(z_{\lambda})-\lambda\frac{ \|z_{\lambda}\|^2}{2}+\lambda\frac{\|x\|^2}{2}=\langle z'_{\lambda},p \rangle-(\tilde f_x)^{**}(z'_{ \lambda})-\lambda\frac{ \|z'_{\lambda}\|^2}{2}+\lambda\frac{\|x\|^2}{2},
\end{eqnarray}
from which we obtain
\begin{equation}\label{omeg}z_{\lambda}=z'_{\lambda}=\nabla_2 L_{\lambda}(x,p) \in \bar \Omega.\end{equation} Therefore, it follows
from (\ref{ine}) that
\begin{eqnarray*}
\lim_{r\rightarrow 0} \frac {L_{r, \lambda}(x,p)-L_{\lambda}(x,p)}{r}=H(\nabla_2
L_{\lambda}(x,p),x).
\end{eqnarray*}
\hfill$\square$ \\
\textbf{End of the proof of Theorem \ref{main.2}:}  For each $\lambda >0, $ $x \in \bar \Omega$ and $p \in \R^N$, we define $S_{\lambda} (x,p)=\nabla_2 L_{\lambda}(x,p).$ It is easy to see that $S_{\lambda} (x,p) \to S_0(x,p)$ where $S_0(x,p)$ is the unique  element with minimal norm in $\partial_2 L_{\infty}(x,p)$ (see Proposition 1.3 in \cite{Barb}).  Set $S(x)= S_0(x, u(x)).$ For each $r>0,$ $\lambda \in [0,1]$ and $x \in \bar \Omega,$ define
\[\eta_r(\lambda, x)= \frac{L_{r,\lambda}(x, u(x)) -L_{\lambda}(x, u(x))}{r}.\]
Note that the function $r \to L_{r,\lambda}(x, u(x))$ is a convex function because it is supremum of a family of linear functions. Thus, for fixed $(x, \lambda)\in \Omega \times [0,1]$,  the function $r \to \eta_r(\lambda, x)$ is non-decreasing. Setting $\eta_0(\lambda, x)$ to be $H(S_\lambda(x), x)$ for $\lambda >0$ and $\eta_0(0, x)=H(S(x), x),$ we have that both  functions $\lambda \to \eta_r(\lambda, x)$ and $\lambda \to \eta_0(\lambda, x)$ are continuous. It follows from Dini's Theorem, that for a fixed $x,$ $\eta_r(\lambda, x)$ converges uniformly to $\eta_0(\lambda, x)$ as $r \to 0$ with respect to $\lambda \in [0,1].$

Note also that thanks to (\ref{omeg}) we have that  $S_{\lambda}, S:\bar \Omega \to \bar \Omega$.  We shall now show that $\int_{\Omega} H(S (x), x) \, dx =0$ for every $H\in {\cal H}$, meaning that $S$ is indeed a self-dual point transformation.  Indeed, by Fatou's lemma we have

\[\lim_{\lambda \to 0} \int_{\Omega} H(S_{\lambda} (x), x) \, dx =\int_{\Omega} H(S (x), x) \, dx.\]
It follows from (\ref{ine}) that
\[\Big |\frac {L_{r, \lambda}(x,p)-L_{\lambda}(x,p)}{r} \Big| \leq \|H\|_{L^\infty(B_R \times B_R)}+|r|. \]
It follows that
\begin{eqnarray*}
 \int_{\Omega} H(S (x), x) \, dx&=&\int_{\Omega} \lim_{\lambda \to 0}  \lim_{r\rightarrow 0^+} \frac {L_{r, \lambda}(x,u(x))-L_{\lambda}(x,u(x))}{r}\, dx\\
&=&\int_{\Omega} \lim_{\lambda \to 0}  \lim_{r\rightarrow 0^+} \eta_r(\lambda, x)\, dx\\
&=&\int_{\Omega}  \lim_{r\rightarrow 0^+} \lim_{\lambda \to 0}  \eta_r(\lambda, x)\, dx \qquad \quad \text{ (due to the uniform convergence) }\\
&=&\int_{\Omega}  \lim_{r\rightarrow 0^+}   \eta_r(0, x)\, dx \\
&=&\lim_{r\rightarrow 0^+} \int_{\Omega}    \eta_r(0, x)\, dx  \qquad \quad \text{ (due to the monotone convergence theorem) } \\
&=&  \lim_{r\rightarrow 0^+} \int_{\Omega} \frac {L_{r}(x,u(x))-L_{\infty}(x,u(x))}{r}\, dx\\
&\geq & 0,
\end{eqnarray*}
from which we have $ \int_{\Omega} H(S (x), x) \, dx \geq 0$.

By the same argument considering $r \to 0^-$, one has $ \int_{\Omega} H(S (x), x) \, dx \leq 0$ and therefore  the latter is indeed zero as desired.

It follows from  the fact that $S(x) \in \partial_2 L_{\infty}(x, u(x))$ together with  $(\tilde f_x)^{**}$ being  the Fenchel dual of $L$ with respect to the second variable (in view of Lemma \ref{wei}) that  $u(x) \in \partial (\tilde f_x)^{**} (S(x)).$ By considering Theorem \ref{diff} in the Appendix, assume that $\Omega'$ is a dense subset of $B_R $ such that  for each $z \in \Omega',$  $\nabla_1 H_{\infty}(z,x)$ exists for all $x \in \bar \Omega.$ Define $\Omega_0=S^{-1}(\Omega'\setminus \partial \Omega)$. Since $S$ is measure preserving we have that $\Omega_0$ is dense in $\bar \Omega.$ We also have for each $x \in \Omega_0$, that $\nabla_1H_{\infty}(S(x), x)$ exists. Since $(\tilde f_x)^{**}(.)= H_{\infty} (., x)$ on $\bar \Omega$ we obtain
\[u(x)=\nabla_{1}H_{\infty}(S(x),x), \qquad \quad \text{ for all } x \in \Omega_0.\]

To complete the proof of Theorem \ref{main.2}, it remains to show that $P_{\infty}=D_{\infty}.$ We already know that $P_{\infty}\geq D_{\infty}.$ To prove the equality it suffices to notice the following:
\begin{eqnarray*}
P_{\infty}= \int_{\Omega} L_{\infty}(x, u(x)) \, dx &=&\int_{\Omega} L_{\infty}(x, u(x)) \, dx+ \int_{\Omega}H_{\infty}(S(x),x) \, dx\\&=&\int_{\Omega} L_{\infty}(x, u(x)) \, dx+ \int_{\Omega}(\tilde f_x)^{**}(S(x)) \, dx
\\&=& \int_{\Omega} \langle u(x), S (x) \rangle \, dx\leq D_\infty.
\end{eqnarray*}

\subsection{Remarks on the uniqueness of the decomposition}

We have seen in example (14) that one cannot expect uniqueness of the involution $S$ in the above decomposition of a given vector field $u$.
We now complete the part of the proof of Theorem \ref{main}, which gives the uniqueness of the involution. (\ref{unique}) on $u$.\\

{\bf Proof of Theorem \ref{main}:} Assume first that $u(x)=\nabla_1H_1(S_1x,x),$ for some Hamiltonian $H_1$ and some selfdual transformation  $S_1$. We shall show that $(H_1, S_1)$ is an ``extremal pair" (i.e., where $D_\infty$ and $S_\infty$ are attained), and that $u(x)=\nabla_1 H_\infty (S_1x, x) $, where $H_\infty$ is the optimal Hamiltonian constructed above. Indeed, let $L$ be the Fenchel-Legendre dual of $H_1$ with respect to the first variable. We have that $L_{H_1}\leq L$ on $\R^N \times \Omega$. It follows that
\[\langle u(x), S_1(x) \rangle \leq L_{H_1}( S_1(x), u(x))+ H_1(S_1 x,x) \leq L(x, u(x))+ H_1(S_1 x,x)=\langle u(x), S(x) \rangle \]
from which we have
\[\langle u(x), S_1(x) \rangle = L_{H_1}(x,u(x))+ H_1(S_1 x,x),\]
and \[\int_\Omega \langle u(x), S_1(x) \rangle \, dx= \int_\Omega L_{H_1}(x, u(x)) \, dx.\]
On the other hand we have
\[\int_\Omega \langle u(x), S_1(x) \rangle \, dx \leq D_\infty=P_\infty \leq \int_\Omega L_{H_1}(x,u(x)) \, dx,\]
which yields
\[\int_\Omega \langle u(x), S_1(x) \rangle \, dx = D_\infty=P_\infty = \int_\Omega L_{H_1}(x,u(x)) \, dx.\]

Now we can  show that $u(x)=\nabla_1 H_\infty (S_1x, x).$ In fact,
\[\int_\Omega \langle u(x), S_1(x) \rangle \, dx=\int_\Omega L_{H_1}(x, u(x)) \, dx=P_\infty=\int_\Omega L_\infty(x,u(x)) \, dx +\int_\Omega H_\infty(S_1(x), x) \, dx,\]
which implies that
\[
\langle u(x), S_1(x) \rangle =L_\infty(x, u(x)) +H_\infty(S_1(x), x)
\]
a.e. on $\Omega$, and hence the desired result.\\

Assume now that the function
$
x\to\langle u(x), y_1-y_2\rangle +\langle u(y_1)-u(y_2), x\rangle
$
has no critical point unless when $y_1=y_2$. Suppose $S_1$, $S_2$ are two transformations such that for $i=1,2$, we have
\begin{equation}
 \big (u(x), - u(S_ix)\big)=\nabla_1H_i(S_ix,x).
\end{equation}
We shall show that $S_1=S_2$ a.e. on $\Omega$. Note first that  the previous argument gives that
\begin{equation}
 \big (u(x), - u(S_ix)\big)=\nabla_1H_\infty (S_ix,x).
\end{equation}
Note also that the function $x\to L_\infty(x,u(x))$ is locally Lipchitz and therefore is differentiable on a subset $\Omega_0$ of full measure. We now show that $S_1=S_2$ on $\Omega_0$.

Indeed, for any $x\in \Omega_0$, $h=0$ is a minimum for the function
\[
h\to L_\infty(x+h, u(x+h)) +H_\infty(S_i(x), x+h)-\langle u(x+h), S_i(x)\rangle.
\]
This implies that
\[
\nabla_2H_\infty (S_1(x), x) -D u(x)S_1(x)=-\frac{d}{dh}L_\infty(x+h,u(x+h))_{h=0}=\nabla_2H_\infty(S_2(x), x)-D u(x)S_2(x),
\]
from which it follows that
\[
D u(x) (S_2(x)-S_1(x))=u(S_1(x))-u(S_2(x)).
\]
The hypothesis then implies that $S_1(x)=S_2(x)$, and $S$ is therefore unique.\\

\begin{remark} Note that the function
\[
h\to L_\infty(x+h,u(x)) +H_\infty(S(x), x+h)-\langle u(x), S(x)\rangle.
\]
has a minimum at $h=0$, from which we have
\[
u(S(x))=-\nabla_2H_\infty(S(x), x) \in \partial_1L_\infty (x,u(x)).
\]
Since, on the other hand, we have $S(x)\in \partial_2L_\infty(x,u(x))$, one obtains
\[
\big ( u(S(x)), S(x)\big)\in \partial L_\infty (x,u(x)).
\]
Now under the hypothesis (\ref{unique}) that ensures uniqueness, the above inclusion becomes
\[
\big( u(S(x)), S(x)\big )= \nabla L_\infty (x,u(x)) \quad {\rm a.e.}\,  x\in \Omega.
\]

\end{remark}

Suppose now that $u$ is a monotone map. We shall show that it satisfies condition (\ref{unique}), This will then yield (3) of Theorem \ref{main}, since $u$ is then a.e., differentiable and the decomposition holds with $S$ being the identity, according to the theorem of Krause.

Indeed, any critical point $\bar x$ of the function
$
x\to\langle u(x), y_1-y_2\rangle +\langle u(y_1)-u(y_2), x\rangle
$
satisfies
\[
D u(\bar x) (y_1-y_2) +u(y_1)-u(y_2)=0,
\]
hence
\[
\langle D u(\bar x) (y_1-y_2), y_1-y_2\rangle +\langle u(y_1)-u(y_2), y_1- y_2\rangle=0.
\]
Since both terms are non-negative, they are equal to zero. If $u$ is strictly monotone,  this cannot happen unless $y_1=y_2$.

\section{Appendix}
\begin{theorem} \label{diff} Let $H$ be a skew-symmetric finite convex-concave function on $\R^{N} \times \R^N $ such that for some $\Lambda >0$, it satisfies
\begin{equation}\label{eee}|H(x_1, y_1)-H(x_2, y_2)| \leq \Lambda \|x_1-x_2\|+\Lambda \|y_1-y_2\| \quad \quad \text{  for all  } (x_1,y_1), (x_2,y_2) \in \R^{N} \times \R^N.\end{equation}
   Let $A \subset \R^N$ be a closed ball and let  $ B \subset \R^N$  be  a compact  subset with non-empty interior. Then, there exists a dense subset $A'$ of $A$ such that for each $ x \in A'$, $\nabla_1 H(x,y)$ exists for all $y \in B.$
\end{theorem}
This result is actually a particular case of a more general result established in \cite{Mo},  where the same conclusion is established  for finite  convex-concave functions on $\R^{n} \times \R^m $ with $n \not=m$ and without condition (\ref{eee}). For the special case $n=m=N$ considered in Theorem  \ref{diff}, the proof can be shortened and we shall provide here a sketch for the reader's convenience. We shall need a few preliminary results.\\
The following definition and theorem can be found in \cite{Bart}.
\begin{definition}\label{def-ar} A sequence $\{f_n\}$ of (scalar-valued) functions on an arbitrary set $X$ is said to converge to $f$ quasi-uniformly on $X,$ if $\{f_n\}$ converges pointwise to $f$ and if, for every $\epsilon >0$ and $L \in \mathbb{N}$, there exists a finite number of indices $n_1, n_2, ..., n_k \geq L$,  such that for each $x \in X$,  at least one of the following inequalities holds:
\[|f_{n_i}(x)-f(x)| < \epsilon, \qquad \quad i=1, 2, ..., k.\]
\end{definition}
\begin{theorem}\label{arzela} If a sequence of functions on a topological space $X$ converges to a continuous limit, then the convergence is quasi-uniform on every compact subset of $X.$ Conversely, if the sequence converges quasi-uniformly on a subset of $X$, then the limit is continuous on this subset.
\end{theorem}
For  $(x,y) \in \R^{N} \times \R^N$, the one sided directional derivative of $H$ at $(x,y)$ with respect to $(u,v)$ is defined as the  limit
\[\nabla H(x,y) (u,v)=\lim_{\lambda \to 0^+} \frac{H(x+\lambda u, y+\lambda v)-H(x,y)}{\lambda}\]
provided such a limit exists. It is standard that the directional derivatives
\[\nabla_1 H(x,y) (u)=\lim_{\lambda \to 0^+} \frac{H(x+\lambda u, y)-H(x,y)}{\lambda}\]
and
\[\nabla_2 H(x,y) (v)=\lim_{\lambda \to 0^+} \frac{H(x, y+\lambda v)-H(x,y)}{\lambda}\]
exist.
The following result is due to T. Rockafellar \cite{Rock}.
\begin{theorem}\label{rock} Let $H$ be a convex-concave function on $\R^{N} \times \R^N. $ Let $C \times D$ be an open convex set on which $H$ is finite. Then for each $(x,y) \in C \times D$, $\nabla H(x,y) (u,v)$ exists and is a finite positively homogeneous convex-concave function of $(u,v)$  on $\R^{N} \times \R^N.$ In fact,
\[\nabla H(x,y) (u,v)=\nabla_1 H(x,y) (u)+\nabla_2 H(x,y) (v).\]
\end{theorem}
For each $\lambda>0,$ define the following functions on $\R^{N} \times \R^N.$
\[H_{\lambda} (u,v)=\frac{H(x+\lambda u, y+\lambda v)-H(x,y)}{\lambda},\]
\[\tilde H_{\lambda} (u,v)=\frac{H(x+\lambda u, y+\lambda v)-H(x,y+\lambda v)}{\lambda},\]
$ H^1_{\lambda} (u)=H_{\lambda} (u,0)$ and $ H^2_{\lambda} (v)=H_{\lambda} (0,v).$ Note that $H^1_{\lambda}$ and $H^2_{\lambda}$ are monotone, so that  by Dini's Theorem,
both  $H^1_{\lambda}(u)$ and $H^2_{\lambda}(v)$ converge uniformly on compact subsets of $\R^{N}$ to $\nabla_1 H(x,y) (u)$ and $\nabla_2 H(x,y) (v)$ respectively. We  have the following properties for $H_{\lambda}$.
\begin{proposition} \label{quasi}The following statements hold:
\begin{enumerate}
\item  $H_{\lambda} (u,v)$ converges uniformly to $\nabla_1 H(x,y) (u)+\nabla_2 H(x,y) (v)$ on compact subsets $A \times B$ of   $\R^{N} \times \R^N.$\\
\item  If for some $u \in \R^{N}$ we have $\nabla_{1}H (x,y)u=-\nabla_1 H(x,y)(-u)$ then for each $v \in B$,
\[
\hbox{$\lim_{\lambda \to 0^+}\nabla_1 H(x,y+\lambda v)u= \nabla_1 H(x,y)u$ uniformly on $ B$.}
\]
\item  $\tilde H_{\lambda} (u,v)$
converges uniformly to $\nabla_1 H(x,y) (u)$ on compact subsets $A \times B$ of   $\R^{N} \times \R^N.$
\end{enumerate}
\end{proposition}
\textbf{Proof.} We first show that for each $\epsilon >0$, there exists $\lambda_0>0$, such that for all $0<\lambda <\lambda_0$, we have
\[ H_{\lambda} (u,v) <\nabla_1 H(x,y) (u)+\nabla_2 H(x,y) (v)+\epsilon, \quad \text{ for all } (u,v) \in A \times B. \]
Then by a dual argument we have
\[ H_{\lambda} (u,v) >\nabla_1 H(x,y) (u)+\nabla_2 H(x,y) (v)-\epsilon, \quad \text{ for all } (u,v) \in A \times B, \]
from which we obtain the desired result in part (1). The difference quotient in the function $H_{\lambda}$ can be expressed as
\[\frac{H(x, y+\lambda v)-H(x,y)}{\lambda}+\frac{H(x+\lambda u, y+\lambda v)-H(x,y+\lambda v)}{\lambda},\]
where the first quotient converges uniformly to $\nabla_2 H(x,y) (v)$ on $B.$ Since $H^1_{\lambda} (u)$ converges uniformly to $\nabla_1 H(x,y) (u)$ on $A,$ there exists $ \alpha>0$ such that
 \[\frac{H(x+\alpha u, y)-H(x,y)}{\alpha}< \nabla_1 H(x,y) (u)+\epsilon.\]
Since $H$ is Lipschitz  with Lipschitz constant $\Lambda>0$, for every $v \in B$, we have
  \[\frac{H(x+\alpha u, y+\lambda v)-H(x,y+\lambda v)}{\alpha }< \nabla_1 H(x,y) (u)+\frac{\epsilon}{2}+\frac{2\lambda \Lambda\|v\|}{\alpha}.\]
Let $\lambda_0$ be small enough such that $\frac{2\lambda_0 \Lambda \|v\|}{\alpha} < \epsilon/2$ for all $v \in B$. For each $0<\lambda < \min\{\lambda_0, \alpha\}$ we have
\begin{eqnarray}\label{rock1}
 \nabla_1 H(x,y) (u)+\epsilon &>& \frac{H(x+\alpha u, y+\lambda v)-H(x,y+\lambda v)}{\alpha }\nonumber\\
&\geq& \frac{H(x+\lambda u, y+\lambda v)-H(x,y+\lambda v)}{\lambda },
\end{eqnarray}
from which part (1) follows.\\
We know prove part (2). Note first that
\[
\frac{H(x+\lambda u, y+\lambda v)-H(x,y+\lambda v)}{\lambda }\geq \nabla_1 H(x,y+\lambda v)u,
\]
 from which together with  (\ref{rock1}),  we have
\begin{eqnarray}\label{rock2}
 \nabla_1 H(x,y) (u)+\epsilon >  \nabla_1 H(x,y+\lambda v)u,
\end{eqnarray}
for every  $0<\lambda < \min \{\lambda_0, \alpha\}$ and  $v \in B.$ By a similar  argument we have
\begin{eqnarray}\label{rock3}
 -\nabla_1 H(x,y) (-u)-\epsilon <  -\nabla_1 H(x,y+\lambda v)(-u),
\end{eqnarray}
It follows from (\ref{rock2}) and (\ref{rock3}) that
\[ -\nabla_1 H(x,y) (-u)-\epsilon< -\nabla_1 H(x,y+\lambda v)(-u) \leq \nabla_1 H(x,y+\lambda v)u < \nabla_1 H(x,y) (u)+\epsilon, \]
and the result follows due to assumption $\nabla_{1}H (x,y)u=-\nabla_1 H(x,y)(-u).$
\\
Part (3) follows from the fact that $\tilde H_{\lambda}(u,v)=H_{\lambda}(u,v)-H^2_{\lambda} (v).$
\hfill $\square$
\begin{proposition}\label{diff33} Fix $u \in \R^{N}.$ There exists  $A' \subset A$ with $A'$ dense in $A$, such that for each $x \in A',$ there exists a dense subset of $B$, say $B_{x,u},$ such that for all $y \in B_{x,u}$ we have  $\nabla_{1}H (x,y)u=-\nabla_1 H(x,y)(-u)$.
\end{proposition}
\textbf{Proof.} Define $F: A \to \R $ by $F(x)=\int_B H(x, z) \, dz.$ Note that $F$ is convex and therefore there exists a dense subset  $A' \subset A$ on which  $F $ is differentiable. For every $x \in A'$ we have
\[F'(x)u=\lim_{\lambda \to 0^+} \frac{\int_B\big [ H(x +\lambda u, z)-H(x, z) \big ] \, dz}{\lambda}=\lim_{\lambda \to 0^-} \frac{\int_B\big [ H(x +\lambda u, z)-H(x, z)\big] \, dz}{\lambda}.\]
due to Lebesgue monotone convergence theorem we have
\[F'(x)u=\int_B\lim_{\lambda \to 0^+} \frac{ H(x +\lambda u, z)-H(x, z)}{\lambda} \, dz=\int_B \lim_{\lambda \to 0^-} \frac{H(x +\lambda u, z)-H(x, z) }{\lambda}\, dz, \]
from which we have
\[\int_B \big [\nabla_1 H(x, z)u +\nabla_1 H(x,z) (-u)\big ] \, dz=0.\]
Since the integrand is nonnegative there exists a dense subset  $B_{x,u}$ of $B$ such that
\[\nabla_1 H(x, z)u +\nabla_1 H(x,z) (-u)=0, \qquad \quad \text{for all } z \in B_{x,u}.\]
\hfill $\square$

\textbf{Proof of Theorem \ref{diff}.} Let $A'$ be as in the above Proposition. Fix $x \in A'.$ We shall show that for all $u \in \R^{N}$ and $y \in B,$ we have $\nabla_1 H(x, y)u +\nabla_1 H(x,y) (-u)=0$ from  which we obtain  $\nabla_1 H(x,y)$ exists for all $y \in B.$
Fix $u \in \R^{N}$ and define $f(y)= \nabla_1 H(x, y)u.$ We first show that $f$ is continuous on $B.$ Note first that $f(y)=\lim_{n \to \infty } f_n (y)$ where
\[f_n (y)=\frac{H(x+ \lambda_n u , y)-H(x,y)}{\lambda_n},\]
and $\lambda_n=1/n.$
We shall show that $f_n$ converges quasi-uniformly to $f$ on $B.$ Fix $\epsilon >0$ and $L \in \mathbb{N}.$ It follows from Proposition \ref{diff33} that there exists a dense subset $B_{x,u}$ of $B$ such  that \[\nabla_1 H(x, y)u +\nabla_1 H(x,y) (-u)=0, \qquad \text{ for all } y \in B_{x,u}.\]

For each  $y \in B_{x,u},$ it follows from Proposition  \ref{quasi} that there exists $n_y> L$ such that
\[\Big |\frac{H(x+\lambda_{n_y} u, y+\lambda_{n_y} v)-H(x,y+\lambda_{n_y} v)}{\lambda_{n_y} }- f(y)\Big | < \frac{\epsilon}{2},\]
and
\[ \big | \nabla_1 H(x,y+\lambda_{n_y} v)u -f(y) \big | < \frac{\epsilon}{2},\]
for every  $ v \in B.$ This implies that
\begin{equation}\label{quazi1}
|f_{n_y}(y +\lambda_{n_y} v)-f(y+\lambda_{n_y} v) |<\epsilon,
\end{equation}
for all $v \in B$. Define $U_y=\{ y +\lambda_{n_y}v; v \in  B\}.$ Since $B_{x,u}$ is dense in $B$ we have \[B \subset \cup_{y \in B_{x,u}} int \big (U_{y}\big).\] $B$ is compact and therefore there exist $y_1, y_2,..., y_k \in B_{x,u}$ such that $B \subset \cup^k_{i=1} int \big(U_{y_i}).$  This together with (\ref{quazi1}) implies that $f_n$ converges to $f$ quasi-uniformly on $B$ and therefore $f$ is continuous.\\

Since $f$ ic continuous and $\nabla_1 H(x, y)u +\nabla_1 H(x,y) (-u)=0$  for almost all $y \in B,$ we indeed have
\[\nabla_1 H(x, y)u +\nabla_1 H(x,y) (-u)=0,\]
for all $y \in B.$ This completes the proof. \hfill $\square$

\end{document}